\documentclass[10pt]{article}
\usepackage{amsfonts}
\usepackage{amsmath}
\usepackage{mathrsfs}
\usepackage{mathrsfs,amscd,amssymb,color,amsthm,amsmath,bm,graphicx,psfrag,subfigure}
\input{epsf}
\makeatletter

\renewcommand{\@seccntformat}[1]{{\csname the#1\endcsname}{\normalsize .}\hspace{.5em}}
\makeatother

\usepackage{indentfirst}

\def \[{\begin{equation}}
\def \]{\end{equation}}

\newtheorem{thm}{Theorem}[section]
\newtheorem{prop}[thm]{Property}

\newtheorem{lem}[thm]{Lemma}
\newtheorem{cor}[thm]{Corollary}

\newtheorem{fac}{Fact}

\newenvironment{wst}
{\setlength{\leftmargini}{1.5\parindent}
 \begin{itemize}
 \setlength{\itemsep}{-1.1mm}}
{\end{itemize}}

\begin{document}

\setlength{\baselineskip}{17pt}
\begin{center}{\Large \bf Relation between the $H$-rank of a mixed graph and the rank\\[5pt] of its underlying graph\footnote{Financially supported by the National Natural Science Foundation of China (Grant Nos. 11671164, 11271149).}}
\vspace{4mm}

{\large Chen Chen$^a$,\ Shuchao Li $^{a,}$\footnote{Corresponding author.\\
\hspace*{5mm}E-mail: chenmath2017@163.com (C. Chen),
\ lscmath@mail.ccnu.edu.cn (S.C.
Li)},\ Minjie Zhang $^b$}\vspace{2mm}

$^a$ Faculty of Mathematics and Statistics,  Central China Normal
University, Wuhan 430079, PR China

$^b$ School of Mathematics and Statistics, Hubei University of Arts and Science, Xiangyang 441053, PR China
\end{center}


\noindent {\bf Abstract}:\ Given a simple graph $G=(V_G, E_G)$ with vertex set $V_G$ and edge set $E_G$, the mixed graph $\widetilde{G}$ is obtained from $G$ by orienting some of its edges. Let $H(\widetilde{G})$ denote the Hermitian adjacency matrix of $\widetilde{G}$ and $A(G)$ be the adjacency matrix of $G$.  The $H$-rank (resp. rank) of $\widetilde{G}$ (resp. $G$), written as $rk(\widetilde{G})$ (resp. $r(G)$), is the rank of $H(\widetilde{G})$ (resp. $A(G)$). Denote by $d(G)$ the dimension of cycle spaces of $G$, that is $d(G) = |E_G|-|V_G|+\omega(G)$, where $\omega(G),$ denotes the number of connected components of $G$. In this paper, we concentrate on the relation between the $H$-rank of $\widetilde{G}$ and the rank of $G$.  We first show that $-2d(G)\leqslant rk(\widetilde{G})-r(G)\leqslant 2d(G)$ for every mixed graph $\widetilde{G}$. Then we characterize all the mixed graphs that attain the above lower (resp. upper) bound. By these obtained results in the current paper, all the main results obtained in \cite{004,1} may be deduced consequently.

\vspace{2mm} \noindent{\bf Keywords}: Mixed graph; $H$-Rank; Lower-optimal; Upper-optimal; $\delta$-transformation

\setcounter{section}{0}
\section{\normalsize Introduction}\setcounter{equation}{0}
In this paper we consider only graphs without loops and multiple edges. Let $G=(V_G,E_G)$ be a graph, where $V_G$ is the vertex set and $E_G$ is the edge set. We call $n:=|V_G|$ the \textit{order} of $G$ and $|E_G|$ the \textit{size} of $G$. A \textit{mixed graph} is a graph where both directed and undirected edges may exist. Thus, a mixed graph, $\widetilde{G},$ is obtained from $G,$ by orienting some of its edges, and we call $G$ the \textit{underlying graph} of $\widetilde{G}$. For convenience, for a mixed graph $\widetilde{G}=(V_{\widetilde{G}}, E_{\widetilde{G}}),$ one has $V_{\widetilde{G}}=V_G$ and the edge set
$E_{\widetilde{G}}$ is the union of $E^0_{\widetilde{G}}$ and $E^1_{\widetilde{G}},$ where $E^0_{\widetilde{G}}$ is the set of undirected edges of $E_{\widetilde{G}}$ and $E^1_{\widetilde{G}}$ is the set of directed edges of $E_{\widetilde{G}}$. We denote an undirected edge by $\{u,v\}$ and a directed edge (or an arc) from $u$ to $v$ by $(u,v)$.

Denote by $P_n, C_n$ and $K_n$ a path, a cycle and a complete graph of order $n$, respectively. 
Let $d(G)$ denote the \textit{dimension} of cycle spaces of a graph $G.$ Then $d(G)=|E_G|-|V_G|+\omega(G)$, here $\omega(G)$ is the number of connected components of $G$. Two distinct edges in a graph $G$ are \textit{independent} if they do not have a common end-vertex in $G$. A set of pairwise independent edges of $G$ is called a \textit{matching}, while a matching with the maximum cardinality is called a \textit{maximum matching}. This maximum cardinality is called the \textit{matching number} of $G$ and written by $m(G).$ A graph is called an \textit{empty graph} if it has no edges.

The \textit{adjacency matrix} $A(G)=(a_{ij})$ of $G$ is an $n\times n$ matrix whose $a_{ij}=1$ if vertices $i$ and $j$ are adjacent and 0 otherwise. The \textit{skew-adjacency matrix} associated to an oriented graph $G^\sigma$, written as $S(G^\sigma)$, is defined to be an $n\times n$ matrix ${(}s_{uv})$ such that $s_{uv} = 1$ if there is an arc from $u$ to $v$, $s_{uv} = -1$ if there is an arc from $v$ to $u$ and $s_{uv} = 0$ otherwise.
The \textit{Hermitian adjacency matrix} of a mixed graph $\widetilde{G}$ is defined to be an $n\times n$ matrix $H(\widetilde{G})=(h_{uv})$ such that $h_{uv}=1$ if $\{u,v\}\in E_{\widetilde{G}}$,  $h_{uv}=i$ if $(u,v) \in E_{\widetilde{G}}$,  $h_{uv}=-i$ if $(v,u)\in E_{\widetilde{G}}$, and $h_{uv}=0$ otherwise. 
This matrix was introduced, independently, by Liu and Li \cite{33} and Guo and Mohar \cite{g1}.
Since $H(\widetilde{G})$ is Hermitian, its eigenvalues are real. The $H$-rank (resp. rank, skew-rank) of $\widetilde{G}$ (resp. $G,\, G^\sigma$), denoted by $rk(\widetilde{G})$ (resp. $r(G), sr(G^\sigma)$), is the rank of $H(\widetilde{G})$ (resp. $A(G), S(G^\sigma)$).

Recently, the study on the $H$-rank and the characteristic polynomial of mixed graphs attracts more and more researchers' attention. Mohar \cite{44} characterized all the mixed graphs with $H$-rank 2 and showed that there are infinitely many mixed graphs with $H$-rank 2 which can not be determined by their $H$-spectrum. Wang et al. \cite{88} identified all the mixed graphs with $H$-rank 3 and showed that all mixed graphs with $H$-rank 3 can be determined by their $H$-spectrum. Liu and Li \cite{33} investigated the properties for characteristic polynomials of mixed graphs and studied the cospectral problems among mixed graphs.  For more properties and applications about the $H$-rank  and eigenvalues of  mixed graphs, we refer the readers to \cite{g5,g1,g3,g4} and the references therein.

Note that, for a mixed graph $\widetilde{G},$ it is possible that $E_{\widetilde{G}}=E^0_{\widetilde{G}}$ or $E_{\widetilde{G}}=E^1_{\widetilde{G}}$. Hence, both oriented graphs and simple graphs can be seen as the special mixed graphs. Wong, Ma and Tian \cite{1} provided a beautiful relation between the skew-rank of an oriented graph and the rank of its underlying graph, which were extended by Huang and Li \cite{004}. Recently, Ma, Wong and Tian \cite{128} determined the relationship between $sr(G^\sigma)$ and the matching number $m(G)$, whereas in \cite{31} they characterized the relationship between rank of $G$ and its number of pendant vertices, from which it may deduce the relationship between the skew rank of $G^\sigma$ and its number of pendant vertices. Huang, Li and Wang \cite{126} established the relationship between $sr(G^\sigma)$ and the independence number of its underlying graph $G$. Very recently, Chen, Huang and Li \cite{025} studied the relation between the $H$-rank of a mixed graph and the matching number of its underlying graph.

Motivated from \cite{025,004,128,1}, it is natural and interesting for us to consider the relation between $rk(\widetilde{G})$ and some other parameters of its underling graph. In this paper, we focus our attention on determining the relation between the $H$-rank of a mixed graph and the rank of its underlying graph. 
\begin{thm}\label{thm1}
Let $\widetilde{G}$ be a mixed graph. Then
\begin{eqnarray}\label{1}
  -2d(G)\leqslant rk(\widetilde{G})-r(G)\leqslant 2d(G).
\end{eqnarray}
\end{thm}

In order to characterize the extremal graphs whose $H$-rank attains upper and lower bounds in (\ref{1}), we first introduce a graph transformation (see also in \cite{21} and \cite{1}).

Let $G$ be a graph with at least one pendant vertex. An operation of deleting a pendant vertex and its adjacent vertex from $G$ is called the $\delta$-\textit{transformation} on $G$. Let $G$ be a graph whose cycles share no common vertices. To switch $G$ to an induced subgraph without pendant vertices, a series of $\delta$-transformation is applied to $G$ as follows: If $G$ has no pendant vertices, then itself is as required and we are done; otherwise, applying a step of $\delta$-transformation on $G$, we obtain an induced subgraph $G_1$ of $G.$ If $G_1$ has no pendant vertices, then $G_1$ is as required and we are done; otherwise, applying a step of $\delta$-transformation on $G_1$, we obtain an induced subgraph $G_2$ of $G$ and so on. We terminate $\delta$-transformations until we obtain an induced subgraph $G_0$ of $G$ that is has no pendant vertices. The resultant subgraph $G_0$ is called a \textit{crucial subgraph} of $G.$

Recall that the \textit{characteristic polynomial} of a mixed graph $\widetilde{G}$ and its underlying graph $G$ is defined, respectively, as
\begin{align}
  \phi(\widetilde{G};\lambda)=&\det\left(\lambda I_n-H(\widetilde{G})\right)=\lambda^n+a_1\lambda^{n-1}+a_2\lambda^{n-2}+\cdots+a_n,\label{eq:2.2} \\
  \psi(G;\lambda)=&\det\left(\lambda I_n-A(G)\right)=\lambda^n+c_1\lambda^{n-1}+c_2\lambda^{n-2}+\cdots+c_n, \label{eq:2.3}
\end{align}
where $I_n$ denotes the identity matrix of order $n$.

Given a mixed cycle $\widetilde{C}$, the \textit{signature} of $\widetilde{C}$, written as $\eta(\widetilde{C})$, is defined as $|f-b|$, where $f$ (resp. $b$) denotes the number of forward (resp. back) oriented edges on $\widetilde{C}$. A subgraph $\widetilde{B}$ of $\widetilde{G}$ is called \textit{basic} if it contains only $\widetilde{K}_2$ or mixed cycles with even signature as components. The \textit{signature} of a basic graph $\widetilde{B}$ is defined as $\eta(\widetilde{B})=\sum \eta(\widetilde{C})$, the sum is over all mixed cycles contained in $\widetilde{B}.$
A subgraph $H$ of $G$ is called an \textit{elementary graph} if each component of $H$ is either an edge or a cycle.
\begin{thm}\label{thm2}
Let $\widetilde{G}$ be a mixed graph. Then $rk(\widetilde{G})-r(G)=2d(G)$ holds if and only if all the following conditions hold for $\widetilde{G}:$
\begin{wst}
\item[{\rm (i)}] the mixed cycles (if any) of $\widetilde{G}$ are pairwise vertex-disjoint;
\item[{\rm (ii)}] for each mixed cycle $\widetilde{C}_l$ of $\widetilde{G},$ $l\equiv 0\pmod {4}$ and either $\eta(\widetilde{C}_l)$ is odd or $\eta(\widetilde{C}_l)\equiv 2\pmod {4};$
\item[{\rm (iii)}] a series of $\delta$-transformations can switch $G,$ the underlying graph of $\widetilde{G},$ to a crucial subgraph, which is the disjoint union of $d(G)$ cycles together with some isolated vertices.
\end{wst}
\end{thm}

\begin{thm}\label{thm3}
Let $\widetilde{G}$ be a mixed graph. Then $rk(\widetilde{G})-r(G)=-2d(G)$ holds if and only if all the following conditions hold for $\widetilde{G}:$
\begin{wst}
\item[{\rm (i)}] the mixed cycles (if any) of $\widetilde{G}$ are pairwise vertex-disjoint;
\item[{\rm (ii)}] for each mixed cycle $\widetilde{C}_l$ of $\widetilde{G}$, $l\equiv 2\pmod {4}$ and $\eta(\widetilde{C}_l)\equiv 2\pmod{4}$;
\item[{\rm (iii)}] a series of $\delta$-transformations can switch $G,$ the underlying graph of $\widetilde{G},$ to a crucial subgraph, which is the disjoint union of $d(G)$ induced cycles together with some isolated vertices.
\end{wst}
\end{thm}

Given a mixed graph $\widetilde{G}=(V_{\widetilde{G}}, E_{\widetilde{G}})$ with $E_{\widetilde{G}}=E^0_{\widetilde{G}}\cup E^1_{\widetilde{G}}.$ If $E^0_{\widetilde{G}}=\emptyset,$ then $\widetilde{G}$ is an oriented graph, $G^\sigma$, for some orientation $\sigma.$ In this case, $H(\widetilde{G})=i S(G^\sigma)$, and therefore $sr(G^\sigma)=rk(\widetilde{G})$.

\setlength{\baselineskip}{17pt}

Let $C_k=x_1 x_2 \ldots x_k x_1$ be a cycle of length $k$.
The \textit{sign} of $C_k^\sigma$ with respect to $\sigma$ is defined to be the sign of $(\prod_{i=1}^{k-1}s_{u_i u_{i+1}}) s_{u_k u_1}$. An even oriented cycle is called \textit{evenly-oriented} (resp. \textit{oddly-oriented}) if its sign is positive (resp. negative).

Together with Theorems \ref{thm1}-\ref{thm3}, we can obtain the following two corollaries, which can be found in \cite{1} and \cite{004}, respectively.
\begin{cor}[\cite{1}]\label{cor5}
Let $G^\sigma$ be a finite oriented graph whose underlying graph is simple. Then $sr(G^\sigma)-r(G)\leqslant 2d(G)$. The equality holds if and only if all the following conditions hold for $G^\sigma:$
\begin{wst}
\item[{\rm (i)}] the cycles (if any) of $G^\sigma$ are pairwise vertex-disjoint;
\item[{\rm (ii)}] each cycle of $G^\sigma$ is oddly-oriented with order a multiple of $4;$
\item[{\rm (iii)}] a series of $\delta$-transformations can switch $G$ to a crucial subgraph of $G,$ which is the disjoint union of $d(G)$ cycles together with some isolated vertices.
\end{wst}
\end{cor}

\begin{cor}[\cite{004}]\label{cor6}
Let $G^\sigma$ be a finite oriented graph whose underlying graph is simple. Then $sr(G^\sigma)-r(G)\geqslant -2d(G)$. The equality holds if and only if all the following conditions hold for $G^\sigma:$
\begin{wst}
\item[{\rm (i)}] the cycles (if any) of $G^\sigma$ are pairwise vertex-disjoint;
\item[{\rm (ii)}] for each cycle $C_l^\sigma$ of $G^\sigma$, $C_l^\sigma$ is evenly-oriented with $l\equiv 2\pmod {4};$
\item[{\rm (iii)}] a series of $\delta$-transformations can switch $G$ to a crucial subgraph of $G,$ which is the disjoint union of $d(G)$ cycles together with some isolated vertices.
\end{wst}
\end{cor}
The remainder of this paper is organized as follows. In Section 2, we list or give some preliminary results which will be used to prove our main results. In Sections 3, we give the proofs of Theorem \ref{thm1}. The proofs of Theorems \ref{thm2} and \ref{thm3} are presented, respectively, in Section 4 and Section 5.
\section{\normalsize Preliminaries}
\setlength{\baselineskip}{17pt}
\setcounter{equation}{0}
Given a mixed graph $\widetilde{G}$, we call $\widetilde{H}$ an \textit{induced subgraph} of $\widetilde{G}$ if $H$ is an
induced subgraph of $G$ and each edge of $\widetilde{H}$ has the same orientation (or non-orientation) as that in $\widetilde{G}$. For $X\subseteq V_{\widetilde{G}},\, \widetilde{G}-X$ is the mixed subgraph obtained from $\widetilde{G}$ by deleting all vertices in $X$ and all incident edges or arcs.  In particular, $\widetilde{G}-\{x\}$ is usually written as $\widetilde{G}-x$ for simplicity. For the sake of clarity, we use the notation $\widetilde{G}-\widetilde{H}$ instead of $\widetilde{G}-V_{\widetilde{H}}$ if $\widetilde{H}$ is an induced subgraph of $\widetilde{G}$.

We call $x$ a \textit{pendant vertex} of $\widetilde{G}$ if it is a vertex of degree one in the underlying graph $G$. Similarly, we call $y$ a \textit{quasi-pendant vertex} of $\widetilde{G}$ if it is adjacent to a vertex of degree one in the underlying graph $G$. For an induced subgraph $\widetilde{K}$ and a vertex $x$ outside $\widetilde{K}$, the induced subgraph of $\widetilde{G}$ with vertex set $V_{\widetilde{K}}\cup\{x\}$ is simply written as $\widetilde{K}+x$. An induced mixed cycle of $\widetilde{G}$ is called a \textit{pendant cycle} if in the underlying graph $G$, this cycle contains a unique vertex of degree 3 and each of its rest vertices is of degree 2.

\subsection{\normalsize Some known lemmas}
In this subsection, we give some known results, which will be used to prove our main results. The first lemma follows immediately from the definitions of the $H$-rank and the second lemma follows by the definition of the matching number.
\begin{lem}[\cite{44}]\label{lem1}
Let $\widetilde{G}$ be a mixed graph. If $\widetilde{G}_1, \widetilde{G}_2, \ldots, \widetilde{G}_k$ are connected components of $\widetilde{G}$, then $rk(\widetilde{G})=\sum\limits_{j=1}^k {rk(\widetilde{G}_j)}.$
\end{lem}

\begin{lem}\label{lem11}
Let $G$ be a simple undirected graph. Then $m(G)-1\leqslant m(G-v)\leqslant m(G)$ for any $v\in V_G$.
\end{lem}

Let $G$ be a graph with pairwise vertex-disjoint cycles and let $T_G$ be obtained from $G$ by contracting each cycle into a vertex, which is called the \textit{contracted vertex}. Then let $[T_G]$ be obtained from $T_G$ by deleting all the contracted vertices and the incident edges (see \cite{1} for details).
\begin{lem}[\cite{22}]\label{lemK}
Let $G$ be a simple undirected graph with at least one cycle. Suppose that all cycles of $G$ are pairwise-disjoint and each cycle is odd, then $m(T_G)=m([T_G])$ if and only if $m(G)=\sum_{C\in \mathcal{C}_G} m(C)+m([T_G]),$ where $\mathcal{C}_G$ denotes the set of all cycles in $G.$
\end{lem}

\begin{lem}[\cite{128}]\label{lem6}
Let $x$ be a pendant vertex of $G$ and $y$ be the neighbor of $x.$ Then $m(G)=m(G-y)+1=m(G-x-y)+1$.
\end{lem}

\begin{lem}[\cite{88}]\label{lem2}
Let $\widetilde{T}$ be a mixed tree. Then $rk(\widetilde{T})=2m(T)$.
\end{lem}
\begin{lem}[\cite{88}]\label{lem3}
Let $x$ be a pendant vertex of $\widetilde{G}$ and $y$ be the neighbor of $x$. Then $rk(\widetilde{G})=rk(\widetilde{G}-x-y)+2$.
\end{lem}
\begin{lem}[\cite{44}]\label{lem4}
Let $\widetilde{G}$ be a mixed graph with $x\in V_{\widetilde{G}}$. Then $rk(\widetilde{G})-2\leqslant rk(\widetilde{G}-x)\leqslant rk(\widetilde{G})$.
\end{lem}
Clearly, Lemmas \ref{lem2}, \ref{lem3} and \ref{lem4} can be easily deduced, respectively, the corresponding results for undirected graphs (see \cite{21,17,1997} for details).

The next result characterizes the relation of the dimension of the cycle spaces between $G$ and a subgraph obtained from $G$ by deleting one vertex.
\begin{lem}[\cite{1}]\label{lem8}
Let $G$ be a graph with $x\in V_G$. Then
\begin{wst}
\item[{\rm (i)}] $d(G)=d(G-x)$ if $x$ lies outside any cycle of $G;$
\item[{\rm (ii)}] $d(G-x)\leqslant d(G)-1$ if $x$ lies on a cycle;
\item[{\rm (iii)}] $d(G-x)\leqslant d(G)-2$ if $x$ is a common vertex of distinct cycles.
\end{wst}
\end{lem}
The following result characterizes the relationship of the ranks between a tree and its subgraph.
\begin{lem}[\cite{128}]\label{lem5}
Let $T$ be a tree with at least one edge.
\begin{wst}
\item[{\rm (i)}] $r(T_1)<r(T),$ where $T_1$ is the subtree obtained from $T$ by deleting all the pendant vertices of $T$.
\item[{\rm (ii)}] If $r(T-W)=r(T)$ for a subset $W$ of $V_T$, then there exist a pendant vertex $v$ of $T$ such that $v$ is not in $W.$
\end{wst}
\end{lem}

\begin{lem}[\cite{33}]\label{lem7}
Let $\mathscr{B}_j$ be the set of basic graphs with $j$ vertices of $\widetilde{G}.$ Then the coefficient $a_j$ defined in $(\ref{eq:2.2})$ is
$$
a_j=\sum_{\widetilde{B}\in\mathscr{B}_j}{(-1)^{\frac{1}{2} \eta(\widetilde{B})+\omega(\widetilde{B})}\cdot2^{c(\widetilde{B})}},\ \ \ \  j=1,2,\ldots,n,
$$
where $\omega(\widetilde{B})$ denotes the number of components of $\widetilde{B}$ and $c(\widetilde{B})$ is the number of cycles in $\widetilde{B}$.
\end{lem}

\begin{lem}[\cite{R}]\label{lemW}
Let $\mathscr{H}_j$ be the set of elementary graphs with $j$ vertices of $G.$ Then the coefficient $c_j$ defined in $(\ref{eq:2.3})$ is
$$
c_j=\sum_{H\in\mathscr{H}_j}{(-1)^{c_1(H)+c(H)}\cdot2^{c(H)}},\ \ \ \  j=1,2,\ldots,n,
$$
where $c_1(H)$ and $c(H)$ denotes the number of components in a subgraph $H$ which are edges and cycles, respectively.
\end{lem}

\begin{lem}[\cite{88}]\label{lem10}
Let $\widetilde{C}_n$ be a mixed cycle with $n$ vertices. Then
$$
rk(\widetilde{C}_n)=\left\{
      \begin{array}{ll}
       n-1 , & \hbox{if $n$ is odd, $\eta(\widetilde{C}_n)$ is odd;} \\[5pt]
       n, & \hbox{if $n$ is odd, $\eta(\widetilde{C}_n)$ is even;} \\[5pt]
       n, & \hbox{if $n$ is even, $\eta(\widetilde{C}_n)$ is odd;} \\[5pt]
       n, & \hbox{if $n$ is even, $n+\eta(\widetilde{C}_n) \equiv 2\pmod{4};$} \\[5pt]
       n-2, & \hbox{if $n$ is even, $n+\eta(\widetilde{C}_n) \equiv 0\pmod{4}.$}
     \end{array}
   \right.
$$
\end{lem}

\begin{lem}[\cite{77}]\label{lem12}
Let $C_n$ be a cycle with $n$ vertices. Then
$r(C_n)=n-2$ if $n\equiv 0\pmod {4},$ and $n$ otherwise.
\end{lem}

\begin{lem}[\cite{22,999,g2}]\label{lem33}
Let $G$ be an undirected graph. Then
\[\label{eq:2.1}
  -2d(G)\leqslant r(G)-2m(G)\leqslant d(G).
\]
The left equality in $(\ref{eq:2.1})$ holds if and only if
all the following conditions hold for $G$:
\begin{wst}
\item[{\rm (i)}] all cycles (if any) of $G$ are pairwise vertex-disjoint;
\item[{\rm (ii)}] the length of each cycle (if any) of $G$ is odd;
\item[{\rm (iii)}] $m(T_G)=m([T_G])$.
\end{wst}
The right equality in $(\ref{eq:2.1})$ holds if and only if
all the following conditions hold for $G$:
\begin{wst}
\item[{\rm (i)}] all cycles (if any) of $G$ are pairwise vertex-disjoint;
\item[{\rm (ii)}] the length of each cycle (if any) of $G$ is a multiple of $4;$
\item[{\rm (iii)}] $m(T_G)=m([T_G])$.
\end{wst}
\end{lem}
\subsection{\normalsize Our preliminaries}
In this subsection, we give some preliminary results, which will be used to prove our main results.
\begin{lem}\label{lem0.2}
Let $\widetilde{G}$ be a mixed graph with $m(T_G)=m([T_G]),$ then each vertex lying on a mixed cycle is not a quasi-pendant vertex of $\widetilde{G}$.
\end{lem}
\begin{proof}
If there exists a quasi-pendant vertex $y$ lying on a mixed cycle of $\widetilde{G}$, let $x$ be the pendant vertex which is adjacent to $y$ and $M$ be a maximum matching of $[T_G]$. Then it is routine to check that $M \cup \{xy\}$ is a matching of $T_G$, thus we get $m(T_G)\geqslant m([T_G])+1$, a contradiction.
\end{proof}

\begin{lem}\label{lem2222}
Let $\widetilde{G}$ be a mixed graph. Then
\begin{eqnarray}\label{2.4}
  -2d(G)\leqslant rk(\widetilde{G})-2m(G)\leqslant d(G).
\end{eqnarray}
Moreover, if the right (resp. left) equality in $(\ref{2.4})$ holds, then every vertex lying on a mixed cycle of $\widetilde{G}$ is not a quasi-pendant vertex.
\end{lem}
\begin{proof}
First we show that $rk(\widetilde{G})\leqslant2m(G)+d(G).$ Let $\phi(\widetilde{G};\lambda)=\sum_{j=0}^na_j\lambda^{n-j}$ be the characteristic polynomial of $\widetilde{G}$. Then it follows from Lemma \ref{lem7} that $a_{j}=0$ for any $j>2m(G)+d(G),$ which is based on the fact that $\widetilde{G}$ contains no basic graphs with $j$ vertices if $j>2m(G)+d(G)$. Consequently, $rk(\widetilde{G})\leqslant 2m(G)+d(G)$. Besides, if $rk(\widetilde{G})=2m(G)+d(G)$ and $u$ is a pendant vertex with neighbor $v.$ Put $\widetilde{G}'=\widetilde{G}-u-v.$ By Lemmas \ref{lem6} and \ref{lem3} one has
\begin{eqnarray}\label{eq:2.00}
m(G')=m(G)-1, \ \ \ rk(\widetilde{G}')=rk(\widetilde{G})-2.
\end{eqnarray}
Combining (\ref{eq:2.00}) with $rk(\widetilde{G})=2m(G)+d(G)$ yields $d(G)=rk(\widetilde{G}')-2m(G').$ On the one hand, by the proof as above, $rk(\widetilde{G}')-2m(G')\leqslant d(G')$, i.e., $d(G)\leqslant d(G').$ On the other hand, $G'$ is a subgraph of $G$, hence $d(G')\leqslant d(G).$ Thus, $d(G)=d(G').$ Consequently, any vertex lying on the mixed cycle of $\widetilde{G}$ is not a quasi-pendant of $\widetilde{G}.$

Now we proceed by induction on $d(G)$ to prove the inequality on the left in (\ref{2.4}). If $d(G)=0$, then $G$ is a forest and the result follows immediately from Lemma~\ref{lem2}. Suppose that $G$ contains at least one cycle, i.e., $d(G) \geqslant 1,$ and let $x$ be a vertex on some cycle of $G$. By Lemma \ref{lem8}, we have
\begin{eqnarray}\label{eq:2.5}
d(G-x)\leqslant d(G)-1.
\end{eqnarray}
By induction hypothesis one has
\begin{eqnarray}\label{eq:2.6}
rk(\widetilde{G}-x)-2m(G-x)\geqslant-2d(G-x).
\end{eqnarray}
According to Lemmas \ref{lem11} and \ref{lem4}, we obtain
\begin{eqnarray}\label{eq:2.7}
m(G-x)\geqslant m(G)-1,\ \ \  rk(\widetilde{G}-x)\leqslant rk(\widetilde{G}).
\end{eqnarray}
Hence, it follows from (\ref{eq:2.5})-(\ref{eq:2.7}) that
\begin{eqnarray*}
rk(\widetilde{G})-2m(G)\geqslant-2d(G).
\end{eqnarray*}
Moreover, if $rk(\widetilde{G})-2m(G)=-2d(G),$ then all inequalities in (\ref{eq:2.5})-(\ref{eq:2.7}) turn into equalities. If $x$ is a quasi-pendant vertex being adjacent to a pendant vertex $z,$ then $rk(\widetilde{G}-x)=rk(\widetilde{G}-z-x)=rk(\widetilde{G})-2<rk(\widetilde{G}),$ a contradiction to $rk(\widetilde{G}-x)= rk(\widetilde{G})$.

This completes the proof.
\end{proof}

\setlength{\baselineskip}{17pt}

\begin{lem}\label{lem555}
Let $\widetilde{G}$ be a {{mixed graph}} containing the unique mixed cycle $\widetilde{C}_l.$ Then $rk(\widetilde{G})=2m(G)+1$ if and only if $l$ is odd, $\eta(\widetilde{C}_l)$ is even and $m(T_G)=m([T_G])$.
\end{lem}
\begin{proof}
Let
\[\label{eq:2.6}
  \phi(\widetilde{G};\lambda)=\sum_{j=0}^na_j\lambda^{n-j}
\] be the characteristic polynomial of $\widetilde{G}$ and $m:=m(G)$. By Lemma \ref{lem2222}, we have $rk(\widetilde{G})\leqslant 2m+1.$ Then combining (\ref{eq:2.6}) with Lemma \ref{lem7} yields
\begin{equation}\label{eq:2.07}
\text{$rk(\widetilde{G})=2m+1$ if and only if $a_{2m+1}\neq 0$.}
\end{equation}

For ``sufficiency", by Lemma \ref{lem7}, we have
$$
a_{2m+1}=\sum_{\widetilde{B}\in\mathscr{B}_{2m+1}}{(-1)^{\frac{1}{2} \eta(\widetilde{B})+\omega({\widetilde{B}})}\cdot2^{c({\widetilde{B}})}}.
$$
Since $l$ is odd and $m(T_G)=m([T_G])$, by Lemma \ref{lemK}, we have $m(G)=m(C_l)+m([T_G])$, i.e.,
$
\frac{l-1}{2}+m([T_G])=m(G),
$
which is equivalent to
\[\label{eq:2.8}
l+2m([T_G])=2m(G)+1.
\]
Let $M_1$ be a maximum matching of $[T_G]$, then $|M_1|=m([T_G]).$ Together with (\ref{eq:2.8}) we obtain that the order of $\widetilde{G}[M_1] \cup \widetilde{C}_l$ is $2m([T_G])+l=2m(G)+1.$ Combining even $\eta(\widetilde{C}_l)$ yields ${\widetilde{M}_1 \cup \widetilde{C}_l} \in \mathscr{B}_{2m+1}$. Then we have
$$
a_{2m+1}=\sum_{\widetilde{B}\in\mathscr{B}_{2m+1}} {(-1)^{\frac{1}{2} \eta(\widetilde{C}_l)+1+\frac{2m+1-l}{2}}\cdot 2^1}=2|\mathscr{B}_{2m+1}| (-1)^{\frac{\eta(\widetilde{C}_l)+3-l}{2}+m}\neq 0.
$$
By Claim 1, we obtain $rk(\widetilde{G})=2m+1.$

For ``necessity", in view of (\ref{eq:2.07}), we have $a_{2m+1}\neq 0.$ Then there exists at least one basic subgraph of order $2m+1.$ Note that each basic subgraph of order $2m+1$ must contain the mixed cycle, $\widetilde{C}_l,$ as its connected component. Thus $l$ is odd and $\eta(\widetilde{C}_l)$ is even.

Next we proceed by induction on $|V_{T_G}|$ to show $m(T_G)=m([T_G]).$ If $|V_{T_G}|=1,$ then $\widetilde{G}\cong \widetilde{C}_l,$ which implies $m(T_G)=m([T_G])=0.$ Now suppose $|V_{T_G}|\geqslant 2.$ If $T_G$ is empty, then $m(T_G)=m([T_G])=0.$ If $T_G$ is non-empty, then there exists a pendant vertex $x$ of $T_G$ which is also a pendant vertex of $\widetilde{G}.$ Let $y$ be the unique neighbor of $x$ and $\widetilde{G}_0=\widetilde{G}-x-y.$ By Lemma \ref{lem2222}, $y\notin V_{C_l}$. Together with Lemmas \ref{lem6} and \ref{lem3} we obtain $rk(\widetilde{G}_0)=rk(\widetilde{G})-2=2m(G)-1=2m(G_0)+1.$ By induction hypothesis, one has
$m(T_{G_0})=m([T_{G_0}]).$ Then it follows from Lemma \ref{lem6} that  $$m(T_G)=m(T_{G_0})+1=m([T_{G_0}])+1=m([T_G]),$$ as desired.
\end{proof}

\begin{lem}\label{lem2.3}
Let $\widetilde{G}$ be a {{mixed graph}} with the unique mixed cycle $\widetilde{C}_l.$ Then $rk(\widetilde{G})=2m(G)-2$ if and only if $l$ is even, $\eta(\widetilde{C}_l)\equiv l\pmod{4}$ and $m(T_G)=m([T_G])$.
\end{lem}
\begin{proof}
For ``sufficiency", we proceed by induction on the order of $T_G$ to prove $rk(\widetilde{G})=2m(G)-2.$ If $|V_{T_G}|=1$, then $\widetilde{G}\cong \widetilde{C}_l.$ Since $l$ is even and $\eta(\widetilde{C}_l)\equiv l\pmod{4}$, by Lemma \ref{lem10}, we have $rk(\widetilde{G})=l-2=2m(G)-2$.

Now suppose that $|V_{T_G}| \geqslant 2$. If $T_G$ is empty, then our result holds immediately. Otherwise, there is a pendant vertex of $T_G$ which is also a pendant vertex of $\widetilde{G}$. Let $x$ be a pendant vertex and $y$ be the unique neighbor of $x$, by Lemma \ref{lem0.2}, $y$ is not on any mixed cycle of $\widetilde{G}$. Denote
$\widetilde{G}_0=\widetilde{G}-x-y$, then it follows from Lemmas \ref{lem6} and \ref{lem3} that
\begin{eqnarray}\label{eq:2.09}
rk(\widetilde{G})=rk(\widetilde{G}_0)+2,\ \ m(G)=m(G_0)+1.
\end{eqnarray}

Note that $|V_{T_{G_0}}|<|V_{T_G}|$ and $\widetilde{G}_0$ is a graph with {{the unique} cycle $\widetilde{C}_l,$} then by induction hypothesis to $\widetilde{G}_0,$ we have
\begin{equation}\label{eq:2.10}
    rk(\widetilde{G}_0)=2m(G_0)-2.
\end{equation}
It follows from (\ref{eq:2.09}) and (\ref{eq:2.10}) that we have $rk(\widetilde{G})=2m({{G}})-2,$ as desired.

For ``necessity", let $\phi(\widetilde{G},\lambda)=\sum_{j=0}^n {a_j \lambda^{n-j}}$ be the characteristic polynomial of $\widetilde{G}$ and $m:=m(G)$. By Lemma~\ref{lem7}, we obtain
\begin{eqnarray}\label{eq:2.11}
a_{2m}&=&\sum_{\widetilde{B}\in \mathscr{B}_{2m}^*} {(-1)^{\frac{1}{2} \eta(\widetilde{B})+\omega(\widetilde{B})} 2^1}+\sum_{\widetilde{M}\in \mathscr{M}} {(-1)^m 2^0}\notag\\
      &=&2\sum_{\widetilde{B}\in \mathscr{B}_{2m}^*} {(-1)^{\frac{\eta(\widetilde{C}_l)-l}{2} +1+m}}+\sum_{\widetilde{M}\in \mathscr{M}} {(-1)^m}\notag\\
      &=& (-1)^m \left(|\mathscr{M}|+2|\mathscr{B}_{2m}^*|(-1)^{\frac{\eta(\widetilde{C}_l)-l}{2}+1}\right),
\end{eqnarray}\setlength{\baselineskip}{17pt}
where $\mathscr{B}_{2m}^*$ denotes the set of all basic subgraphs of order $2m$ containing  $\widetilde{C}_l$ as its connected component and $\mathscr{M}$ is the set of all maximum matchings of $\widetilde{G}$.\setlength{\baselineskip}{17pt}

Since $rk(\widetilde{G})=2m-2$, we have $a_{2m}=0.$ In view of (\ref{eq:2.11}), we have $\frac{\eta(\widetilde{C}_l)-l}{2}$ is even and $|\mathscr{B}_{2m}|=\frac{|\mathscr{M}|}{2}>0$, which implies that $\eta(\widetilde{C}_l)$ is even and $\eta(\widetilde{C}_l)\equiv l\pmod{4}$. Consequently, $l$ is even and $\eta(\widetilde{C}_l)\equiv l \pmod{4}.$

Next we proceed by induction on $|V_{T_G}|$ to show $m(T_G)=m([T_G]).$ If $|V_{T_G}|=1$, then $\widetilde{G}\cong \widetilde{C}_l$ and $m(T_G)=m([T_G])=0$.

Now assume that $|V_{T_G}|\geqslant 2$. {{If $T_G$ is empty, then $m(T_G)=m([T_G])=0.$ If $T_G$ is non-empty,}} then there exists a pendant vertex $u$ of $T_G$ which is also a pendant vertex of $\widetilde{G}$. Let $v$ be the unique neighbor of $u$, then it follows from Lemma \ref{lem2222} that $v$ is not on any mixed cycle of $\widetilde{G}$. Denote $\widetilde{G}_0=\widetilde{G}-u-v$. Then the condition $rk(\widetilde{G})=2m(G)-2$ together with Lemmas \ref{lem6} and \ref{lem3} lead to
$
rk(\widetilde{G}_0)=2m(G_0)-2.
$
By induction hypothesis, one has $m(T_{G_0})=m([T_{G_0}]).$
Therefore, by Lemma \ref{lem6} we have
$$
m(T_G)=m(T_{G_0})+1=m([T_{G_0}])+1=m([T_G]),
$$
as desired.
\end{proof}
\subsection{\normalsize Some basic facts}
A mixed graph $\widetilde{G}$ is called the \textit{upper-optimal} (resp. \textit{lower-optimal}) if $rk(\widetilde{G})-r(G)$ attains the upper bound (resp. lower bound) $2d(G)$ (resp. $-2d(G)$). 
By the definition of the upper-optimal (resp. lower optimal) and Lemma 2.1, the next fact holds immediately.
\setcounter{fac}{0}
\begin{fac}\label{fac1}
$\widetilde{G}$ is upper-optimal (resp. lower-optimal) if and only if all connected components of $\widetilde{G}$ are upper-optimal (resp. lower-optimal).
\end{fac}

By Lemmas \ref{lem10} and \ref{lem12}, the following fact follows immediately.
\begin{fac}\label{fac2}
For a mixed cycle $\widetilde{C}_l$, it is lower-optimal if and only if $l\equiv 2\pmod 4$ and $\eta(\widetilde{C}_l)\equiv 2\pmod 4,$ whereas it is upper-optimal if and only if $l\equiv 0\pmod 4$ and either $\eta(\widetilde{C}_l)$ is odd or $\eta(\widetilde{C}_l)\equiv 2\pmod 4.$
\end{fac}
\section{\normalsize Proof of Theorem 1.1}\setcounter{equation}{0}
\setlength{\baselineskip}{17pt}
In this section, we give the proof for Theorem 1.1.\vspace{2mm}

\noindent{\bf Proof of Theorem 1.1}\ \
We proceed by induction on $d(G)$ to show our result.  If $d(G)=0$, then $\widetilde{G}$ is a mixed forest and the result follows immediately from Lemma \ref{lem2}. Now suppose that $\widetilde{G}$ contains at least one mixed cycle, i.e., $d(G) \geqslant 1,$ and let $x$ be a vertex on a mixed cycle of $\widetilde{G}$. By Lemma \ref{lem8}, we have
\begin{eqnarray}\label{eq:3.01}
d(G-x)\leqslant d(G)-1.
\end{eqnarray}
Applying the induction hypothesis on $\widetilde{G}-x$ yields
\begin{eqnarray}\label{eq:3.02}
rk(\widetilde{G}-x)-r(G-x)\geqslant-2d(G-x)
\end{eqnarray}
and
\begin{eqnarray}\label{eq:3.2}
rk(\widetilde{G}-x)-r(G-x)\leqslant 2d(G-x).
\end{eqnarray}
According to Lemma \ref{lem4}, we obtain
\begin{eqnarray}\label{eq:3.03}
rk(\widetilde{G})\geqslant rk(\widetilde{G}-x),\ \ \  r(G)\leqslant r(G-x)+2
\end{eqnarray}
and
\begin{eqnarray}\label{eq:3.3}
rk(\widetilde{G})\leqslant rk(\widetilde{G}-x)+2,\ \ \  r(G)\geqslant r(G-x).
\end{eqnarray}
Hence, it follows from (\ref{eq:3.01})-(\ref{eq:3.3}) that
\begin{eqnarray*}
rk(\widetilde{G})-r(G)\geqslant-2d(G)
\end{eqnarray*}
and
\begin{eqnarray*}
rk(\widetilde{G})-r(G)\leqslant 2d(G),
\end{eqnarray*}
as desired.  \qed

\section{\normalsize Proof of Theorem 1.2}\setcounter{equation}{0}
Recall that a mixed graph $\widetilde{G}$ is \textit{upper-optimal} if $rk(\widetilde{G})-r(G)$ attains the upper bound $2d(G)$. In this section, we first give some fundamental characterization of upper-optimal mixed graphs. Then we give the proof for Theorem~1.2.
\begin{prop}\label{lem4.1}
Let $x$ be a vertex of $\widetilde{G}$ lying on a mixed cycle. If $\widetilde{G}$ is upper-optimal, then
\begin{wst}
\item[{\rm (i)}] $rk(\widetilde{G})=rk(\widetilde{G}-x)+2;$
\item[{\rm (ii)}] $\widetilde{G}-x$ is upper-optimal;
\item[{\rm (iii)}] $d(G)=d(G-x)+1;$
\item[{\rm (iv)}] $r(G)=r(G-x);$
\item[{\rm (v)}] $x$ lies on just one mixed cycle of $\widetilde{G}$ and $x$ is not a quasi-pendant vertex of $\widetilde{G}$.
\end{wst}
\end{prop}
\begin{proof}
Note that $\widetilde{G}$ is upper-optimal. Together with the proof of (\ref{1}), each of inequalities in (\ref{eq:3.01}), (\ref{eq:3.2}) and (\ref{eq:3.3}) must actually be an equality. Thus, (i)-(iv) hold.

Now we show (v). In fact, by (iii) and Lemma \ref{lem8}(iii), we obtain that $x$ must lie on just one mixed cycle of $\widetilde{G}$. If $x$ is a quasi-pendant vertex being adjacent to a pendant vertex $y$, then by Lemma \ref{lem3}, we have  $r(G-x)=r(G-x-y)=r(G)-2$,  a contradiction to (iv). This completes the proof of (v).
\end{proof}


\begin{prop}\label{lem4.2}
Let $\widetilde{G}$ be a mixed graph containing a pendant vertex $x$ with neighbor $y$. Put $\widetilde{G}'=\widetilde{G}-x-y$. Then $\widetilde{G}$ is upper-optimal if and only if $y$ is not on any mixed cycle of $\widetilde{G}$ and $\widetilde{G}'$ is upper-optimal.
\end{prop}
\begin{proof}
For ``sufficiency", we know that $y$ is not on any cycle of $G$, by Lemma~\ref{lem8},
\begin{eqnarray}\label{eq:4.3}
d(G')=d(G).
\end{eqnarray}
It follows from Lemma \ref{lem3} that
\begin{eqnarray}\label{eq:4.4}
rk(\widetilde{G}')=rk(\widetilde{G})-2, \ \  r(G')=r(G)-2.
\end{eqnarray}
Combing the upper-optimal condition of $\widetilde{G}'$ with (\ref{eq:4.3}) and (\ref{eq:4.4}), we have
$$
rk(\widetilde{G})-r(G)=2d(G),
$$
as desired.

For ``necessity", since $\widetilde{G}$ is upper-optimal,
\begin{eqnarray}\label{eq:4.4s}
rk(\widetilde{G})-r(G)=2d(G).
\end{eqnarray}
Substituting (\ref{eq:4.4}) into  (\ref{eq:4.4s}) yields
$
2d(G)=rk(\widetilde{G}')-r(G').
$
Then in view of Theorem \ref{thm1}, we have $2d(G)=2d(G')=rk(\widetilde{G}')-r(G')$.  Consequently,  $y$ is not on any mixed cycle of $\widetilde{G}$ and $\widetilde{G}'$ is upper-optimal.

This completes the proof.
\end{proof}

\setlength{\baselineskip}{17pt}

\begin{prop}\label{lem4.3}
Let $\widetilde{G}$ be a {{mixed graph}} containing the unique mixed cycle $\widetilde{C}_l.$ Then $\widetilde{G}$ is upper-optimal if and only if all of the following conditions hold for $\widetilde{G}:$
\begin{wst}
\item[{\rm (i)}] $l\equiv 0\pmod 4;$
\item[{\rm (ii)}] either $\eta(\widetilde{C}_l)$ is odd or $\eta(\widetilde{C}_l)\equiv 2\pmod 4;$
\item[{\rm (iii)}] $r(T_G)=r([T_G])$.
\end{wst}
\end{prop}
\begin{proof}
For ``sufficiency", we proceed by induction on $|V_{T_G}|$ to show that $\widetilde{G}$ is upper-optimal. If $|V_{T_G}|=1,$ then $G\cong C_l$, by Fact 3, it is straightforward to check that $rk(\widetilde{G})-r(G)=2.$ Now we assume that $|V_{T_G}|\geqslant 2,$ if $T_G$ is an empty graph, then $\widetilde{G}$ consists of {{the unique mixed cycle $\widetilde{C}_l$}} and isolated vertices. Thus by Facts~1 and 3, $\widetilde{G}$ is upper-optimal.

Now we consider that $T_G$ is non-empty. Together with (iii) and Lemma \ref{lem2}, we obtain $m(T_G)=m([T_G])$. Note that $T_G$ contains a pendant vertex, say $x$. Clearly, it is also a pendant vertex of $\widetilde{G}.$ Let $y$ be the neighbor of $x.$ Then, by Lemma \ref{lem0.2}, $y$ is not on any cycle of $G.$ Put $\widetilde{G}'=\widetilde{G}-x-y.$ Then by Lemma \ref{lem3}, $r(T_{G'})=r(T_G)-2=r([T_G])-2=r([T_{G'}]).$ Applying the induction hypothesis on $\widetilde{G}'$, we have $\widetilde{G}'$ is upper-optimal.
Combining Property \ref{lem4.2} yields that $\widetilde{G}$ is upper-optimal.

For ``necessity", let $m=m(G)$ be the matching number of $G.$ Combining Lemmas \ref{lem33} and \ref{lem2222}, we obtain $2m-2 \leqslant rk(\widetilde{G})\leqslant 2m+1$ and $2m-2 \leqslant r(G)\leqslant 2m+1.$ Note that $\widetilde{G}$ is upper-optimal, i.e. $rk(\widetilde{G})- r(G)=2$. Hence, we proceed by considering the following two possible cases. 

The first case is $rk(\widetilde{G})=2m+1$ and $r(G)=2m-1.$  In this case, by Lemma \ref{lem555}, we have $l$ is odd. Let $\psi (G,\lambda)=\sum\limits_{i=0}^n c_i\lambda^{n-i}$ be the characteristic polynomial of $G.$ Then by Lemma \ref{lemW} we obtain $c_{2m}=\sum\limits_{M\in \mathscr{M}} (-1)^{2m} 2^0=|\mathscr{M}|\neq 0,$ where $\mathscr{M}$ denotes the set of all the maximum matchings of $G.$ Consequently, $r(G)\geqslant 2m,$ a conrtadiction.

The rest case is $rk(\widetilde{G})=2m$ and $r(G)=2m-2.$ In view of Lemma \ref{lem33} we have $l\equiv 0\pmod 4$ and $m(T_G)=m([T_G])$. Together with Lemma \ref{lem2}, we have $r(T_G)=r([T_G]).$ Hence, (i) and (iii) hold. If (ii) does not hold, then $\eta(\widetilde{C}_l)\equiv 0\pmod {4}$. By Lemma \ref{lem2.3}, $rk(\widetilde{G})=2m-2$, a contradiction.

This completes the proof.
\end{proof}

\begin{prop}\label{lem4.4}
Let $\widetilde{G}$ be obtained from mixed graphs $\widetilde{C}_q$ and $\widetilde{H}$ by connecting them with a mixed edge $\tilde{e}=xy$, where $x\in V_{\widetilde{C}_q}$ and $y\in V_{\widetilde{H}}.$ Put $\widetilde{K}:=\widetilde{H}+x.$ If $\widetilde{G}$ is upper-optimal, then
\begin{wst}
\item[{\rm (i)}] for each mixed cycle $\widetilde{C}_l$ of $\widetilde{G},$ $l$ is a multiple of $4$ and either $\eta(\widetilde{C}_l)$ is odd or $\eta(\widetilde{C}_l)\equiv 2\pmod {4}$;
\item[{\rm (ii)}] both $\widetilde{K}$ and $\widetilde{H}$ are upper-optimal;
\item[{\rm (iii)}] $rk(\widetilde{K})=rk(\widetilde{H})$ and $r(K)=r(H);$
\item[{\rm (iv)}] $rk(\widetilde{G})=q+rk(\widetilde{K})$ and $r(G)=q+r(K)-2.$
\end{wst}
\end{prop}
\begin{proof}
(i)\ We begin by induction on $d(G)$ to show our result. If $d(G)=1$, then $\widetilde{G}$ contains a unique mixed cycle and (i) follows directly from Property \ref{lem4.3}. Now suppose that $d(G)\geqslant 2$. Then $\widetilde{H}$ contains at least one mixed cycle. Let $u$ be a vertex lying on some mixed cycle of $\widetilde{H}$. By Property \ref{lem4.1}(ii), $\widetilde{G}_0=\widetilde{G}-u$ is upper-optimal. Since $d(G_0)<d(G)$, by induction each mixed cycle in $\widetilde{G}_0$, including $\widetilde{C}_q$, satisfies (i). By a similar discussion as above, each mixed cycle in $\widetilde{G}-x$ satisfies (i), i.e., all the mixed cycles in $\widetilde{H}$ satisfy (i). That is to say, each mixed cycle in $\widetilde{G}$ satisfies (i).

(ii)\ By (i) we know that $q$ is a multiple of $4$. Then let $C_q=xx_2x_3\dots x_qx$. As $\widetilde{G}$ is upper-optimal, by Property \ref{lem4.1} both $\widetilde{G}-x_2$ and $\widetilde{G}-x$ are upper-optimal. Together with Fact \ref{fac1}, $\widetilde{H}$ is also upper-optimal. Let $\widetilde{G}_1:=\widetilde{G}-\{x_2,x_3,x_4\}$. Note that $x_3$ (resp. $x_4$) is the pendant vertex (resp. quasi-pendant vertex) of $\widetilde{G}-\{x_2\}.$ Then in view of Lemma \ref{lem3}, we have
\begin{eqnarray}\label{eq:4.06}
rk(\widetilde{G}_1)=rk(\widetilde{G}-x_2)-2,\ \ r(G_1)=r(G-x_2)-2.
\end{eqnarray}
As $\widetilde{G}-x_2$ is upper-optimal, we obtain
\begin{eqnarray}\label{eq:4.07}
rk(\widetilde{G}-x_2)-r(G-x_2)=2d(G-x_2).
\end{eqnarray}
Combining (\ref{eq:4.06}) and (\ref{eq:4.07}) with the fact that $d(G_1)=d(G-x_2)$ yields
$
rk(\widetilde{G}_1)-r(G_1)=2d(G_1),
$
i.e.,  $\widetilde{G}_1$  is upper-optimal. Repeating such process, after $\frac{q-2}{2}$ steps, the resultant graph is $\widetilde{G}-\{x_2,x_3,\ldots,x_{q}\}=\widetilde{K},$ which is also upper-optimal.

(iii) and (iv)\ By Lemma \ref{lem3} and Property \ref{lem4.1}, one has
\begin{eqnarray}\label{eq:4.7}
&&rk(\widetilde{G})=rk(\widetilde{G}-x_2)+2=rk(\widetilde{K})+q,\ \ \ \ \
r(G)=r(G-x_2)=r(K)+q-2
\end{eqnarray}
and
\begin{eqnarray}\label{eq:4.9}
&& rk(\widetilde{G})=rk(\widetilde{G}-x)+2=rk(\widetilde{H})+q, \ \ \ \ \ r(G)=r(G-x)=r(H)+q-2.
\end{eqnarray}
Together with (\ref{eq:4.7}) and (\ref{eq:4.9}), (iii) and (iv) hold.

This completes the proof.
\end{proof}

\begin{prop}\label{lem4.6}
Let $\widetilde{G}$ be a mixed graph. If $\widetilde{G}$ is upper-optimal, then
\begin{wst}
\item[{\rm (i)}] the cycles (if any) of $G$ are pairwise vertex-disjoint;
\item[{\rm (ii)}] for each mixed cycle $\widetilde{C}_l$ of $\widetilde{G},$ $l$ is a multiple of $4$ and either $\eta(\widetilde{C}_l)$ is odd or $\eta(\widetilde{C}_l)\equiv 2\pmod {4}$;
\item[{\rm (iii)}] $r(G)=r(T_G)+\sum_{O\in \mathcal{C}_G} (|V_O|-2),$ where $\mathcal{C}_G$ denotes the set of cycles of $G;$
\item[{\rm (iv)}] $r(T_G)=r([T_G]).$
\end{wst}
\end{prop}
\begin{proof}
(i)\ By Property \ref{lem4.1}(v), (i) follows immediately.

(ii)-(iv)\ We proceed by induction on ${{|V_G|}}$ to show (ii)-(iv) simultaneously. If $|V_G|=1,$ then (ii)-(iv) hold trivially. Now suppose $|V_G|\geqslant 2.$ By (i) all cycles (if any) of $G$ are pairwise vertex-disjoint. Hence $T_G$ is well defined. If $T_G$ is empty, then $G$ consists of distinct cycles and isolated vertices. Therefore, by Facts \ref{fac1} and \ref{fac2}, (ii)-(iv) hold immediately. If $T_G$ is non-empty, then $T_G$ has a pendant vertex, say $x.$ If $x$ doesn't lie on any cycle of $G,$ then $x$ is also a pendant vertex of $G.$ If $x$ lies on some cycle of $G,$ then this cycle must be a pendant cycle of $G$. We proceed by considering the following two possible cases to show our result.

{\bf{Case 1.}}\ $x$ is a pendant vertex of $G.$ In this case, let $y$ be the unique neighbor of $x$ and put $\widetilde{G}':=\widetilde{G}-x-y.$ Then by Property \ref{lem4.2}, $y$ is not on any cycle of $G$ and $\widetilde{G}'$ is also upper-optimal. Applying induction hypothesis to $\widetilde{G}'$ yields
\begin{wst}
  \item[{\rm (a)}] for each mixed cycle $\widetilde{C}_l$ of $\widetilde{G}',$ $l$ is a multiple of 4 and either $\eta(\widetilde{C}_l)$ is odd or $\eta(\widetilde{C}_l)\equiv 2\pmod {4}$;
  \item[{\rm (b)}] $r(G')=r(T_{G'})+\sum_{O\in \mathcal{C}_{G'}} (|V_O|-2);$
  \item[{\rm (c)}] $r(T_{G'})=r([T_{G'}]).$
\end{wst}

Note that each of cycles in $\widetilde{G}$ is also that of $\widetilde{G}'.$ Hence, Assertion (a) implies that each cycle of $\widetilde{G}$ satisfies (ii). Clearly, $x$ is a pendant vertex of $T_G$ (resp., $[T_G]$) whose neighbor is $y$. Remind that $T_{G'}=T_G-x-y$ and $[T_{G'}]=[T_G]-x-y$. Hence, together with Lemma \ref{lem3} and (b), we have
$$
r(G)=r(G')+2=[r(T_{G'})+\sum_{O\in \mathcal{C}_{G'}} (|V_O|-2)]+2=r(T_G)+\sum_{O\in \mathcal{C}_G} (|V_O|-2),
$$
i.e., (iii) holds. Combining Lemma \ref{lem3} and (c) yields
$
r(T_G)=r(T_{G'})+2=r([T_{G'}])+2=r([T_G]).
$
Hence, (iv) holds.

{\bf{Case 2.}}\ $x$ lies on a pendant cycle, say $C_q,$ of $G.$ In this case, let $u$ be the unique vertex of $C_q$ of degree 3. Put $\widetilde{H}:=\widetilde{G}-\widetilde{C}_q$ and $\widetilde{K}:=\widetilde{H}+u.$ By Property \ref{lem4.4}, (ii) holds, and both $\widetilde{K}$ and $\widetilde{H}$ are upper-optimal. Furthermore,
\begin{eqnarray}\label{eq:4.11}
r(G)=r(K)+q-2.
\end{eqnarray}
Applying induction to $\widetilde{K}$ yields
\begin{eqnarray}\label{eq:4.12}
r(K)=r(T_K)+\sum_{O\in \mathcal{C}_{K}} (|V_O|-2).
\end{eqnarray}
Note that $T_G\cong T_K$ and $|V_{C_q}|=q.$ Combining (\ref{eq:4.11}) with (\ref{eq:4.12}) yields
$$
r(G)=r(T_G)+\sum_{O\in \mathcal{C}_{G}} (|V_O|-2).
$$
If we apply induction to $\widetilde{H}$, then we see that
\begin{eqnarray}\label{eq:4.13}
r(H)=r(T_H)+\sum_{O\in \mathcal{C}_{H}} (|V_O|-2), \ \ \ r(T_H)=r([T_H]).
\end{eqnarray}
By Property \ref{lem4.4}(iii), one has $r(K)=r(H).$ Combining (\ref{eq:4.12}) with (\ref{eq:4.13}) produces
$r(T_K)=r(T_H).$ Together with $[T_G]\cong [T_H]$, we obtain that $r(T_G)=r(T_K)=r(T_H)=r([T_H])=r([T_G]).$

This completes the proof.
\end{proof}

\noindent{\bf Proof of Theorem 1.2}\ \ Now we come back to characterize the sufficient and necessary conditions for the equality holding on the right of (\ref{1}).

For ``sufficiency", suppose that $\widetilde{G}$ satisfies all the conditions of (i)-(iii). Suppose that we apply, repeatedly, $\delta$-transformations for $t$ times to switch $G$ to a crucial subgraph $G_0,$ which is consisted of $d(G)$ disjoint cycles and some isolated vertices. Note that each $\delta$-transformation decreases the $H$-rank of a mixed graph by 2, so does the rank of its underlying graph. Then
\begin{eqnarray}\label{eq:4.8}
rk(\widetilde{G})=2t+rk(\widetilde{G}_0), \ \ \ r(G)=2t+r(G_0).
\end{eqnarray}
By Lemmas \ref{lem1}, \ref{lem10} and \ref{lem12}, one has
\begin{eqnarray}\label{eq:4.08}
rk(\widetilde{G}_0)=\sum_{\widetilde{O}\in \mathcal{C}_{\widetilde{G}_0}} rk(\widetilde{O})=\sum_{O\in \mathcal{C}_{G_0}} |V_O|=\sum_{O\in \mathcal{C}_{G_0}} (r(O)+2)=r(G_0)+2d(G).
\end{eqnarray}
Substituting (\ref{eq:4.8}) into (\ref{eq:4.08}) yields $rk(\widetilde{G})-r(G)=2d(G),$ as desired.

For ``necessity", (i) and (ii) follow directly by Property \ref{lem4.6}. Hence, in order to complete the proof, it suffices to show (iii). We proceed by induction on $|V_{\widetilde{G}}|$ to prove (iii).

If $|V_{\widetilde{G}}|=1,$ then (iii) holds trivially. Suppose that (iii) holds for every upper-optimal mixed graph $\widetilde{G}$ with $|V_{\widetilde{G}}|<n$. Now let $\widetilde{G}$ be an upper-optimal mixed graph of order $n\geqslant 2.$

If $T_G$ is empty, then $\widetilde{G}$ is the disjoint union of $d(G)$ cycles along with some isolated vertices. Then (iii) holds apparently. If $T_G$ is non-empty, then by Property \ref{lem4.6} we have $r(T_G)=r([T_G])$. Hence, by Lemma \ref{lem5}(ii), $T_G$ contains a pendant vertex, say $x,$ which is also a pendant vertex of $G$. Let $y$ be the neighbor of $x$ and put $\widetilde{G}_1:=\widetilde{G}-x-y.$ By Property \ref{lem4.2}, $y$ is not on any mixed cycle of $\widetilde{G}$ and $\widetilde{G}_1$ is upper-optimal. If we apply induction to $\widetilde{G}_1$ we see that a series of $\delta$-transformations can switch ${G}_1$ to a crucial subgraph $G_0,$ which consists of $d(G_1)$ disjoint cycles together with some isolated vertices. Note that $d(G)=d(G_1)$. Thus, a series of $\delta$-transformations can switch $G$ to the crucial subgraph $G_0$, which is the disjoint union of $d(G)$ cycles together with some isolated vertices.

This completes the proof. \qed
\section{\normalsize Proof of Theorem \ref{thm3}}\setcounter{equation}{0}
Recall that a mixed graph $\widetilde{G}$ is \textit{lower-optimal} if $rk(\widetilde{G})-r(G)$ attains the lower bound $-2d(G)$. In this section, we first give some fundamental characterization of lower-optimal mixed graphs. Then we give the proof for Theorem~1.3.
\begin{prop}\label{lem5.1}
Let $x$ be a vertex of $\widetilde{G}$ lying on a mixed cycle. If $\widetilde{G}$ is lower-optimal, then
\begin{wst}
\item[{\rm (i)}] $rk(\widetilde{G})=rk(\widetilde{G}-x);$
\item[{\rm (ii)}] $\widetilde{G}-x$ is lower-optimal;
\item[{\rm (iii)}] $d(G)=d(G-x)+1;$
\item[{\rm (iv)}] $r(G)=r(G-x)+2;$
\item[{\rm (v)}] $x$ lies on just one mixed cycle of $\widetilde{G}$ and $x$ is not a quasi-pendant vertex of $\widetilde{G}$.
\end{wst}
\end{prop}
\begin{proof}
Note that $\widetilde{G}$ is lower-optimal. Together with the proof of (\ref{1}), each of inequalities in (\ref{eq:3.01}), (\ref{eq:3.02}) and (\ref{eq:3.03}) must actually be an equality. Thus, (i)-(iv) hold.

Now we show (v). In fact, by (iii) and Lemma \ref{lem8}(iii), we obtain that $x$ must lie on just one mixed cycle of $\widetilde{G}$. If $x$ is a quasi-pendant vertex being adjacent to a pendant vertex, say $y$, then by Lemma \ref{lem3}, we have  $rk(\widetilde{G}-x)=rk(\widetilde{G}-x-y)=rk(\widetilde{G})-2$,  a contradiction to (i). This completes the proof of (v).
\end{proof}

\begin{prop}\label{lem5.02}
Let $\widetilde{G}$ be a mixed graph containing a pendant vertex $x$ whose neighbor is $y$. Put $\widetilde{G}'=\widetilde{G}-x-y$. Then $\widetilde{G}$ is lower-optimal if and only if $y$ is not on any mixed cycle of $\widetilde{G}$ and $\widetilde{G}'$ is lower-optimal.
\end{prop}
\begin{proof}
By Lemma \ref{lem3} one has
\begin{eqnarray}\label{eq:5.4}
rk(\widetilde{G}')=rk(\widetilde{G})-2, \ \  r(G')=r(G)-2.
\end{eqnarray}

For the sufficiency, we know that $y$ is not on any cycle of $G$, by Lemma~\ref{lem8},
\begin{eqnarray}\label{eq:5.3}
d(G')=d(G).
\end{eqnarray}
Combing the lower-optimal condition of $\widetilde{G}'$ with (\ref{eq:5.4}) and (\ref{eq:5.3}), we have
$$
rk(\widetilde{G})-r(G)=-2d(G),
$$
as desired.

For ``necessity", since $\widetilde{G}$ is lower-optimal, we get
\begin{eqnarray}\label{eq:5.4s}
rk(\widetilde{G})-r(G)=-2d(G).
\end{eqnarray}
Substituting (\ref{eq:5.4}) into  (\ref{eq:5.4s}) yields
$
-2d(G)=rk(\widetilde{G}')-r(G').
$
Then in view of Theorem \ref{thm1}, we have $-2d(G)=-2d(G')=rk(\widetilde{G}')-r(G')$.  Consequently,  $y$ is not on any mixed cycle of $\widetilde{G}$ and $\widetilde{G}'$ is lower-optimal.

This completes the proof.
\end{proof}

\setlength{\baselineskip}{17pt}

\begin{prop}\label{lem5.3}
Let $\widetilde{G}$ be a {{mixed graph}} containing the unique mixed cycle $\widetilde{C}_l.$ Then $\widetilde{G}$ is lower-optimal if and only if all of the following conditions hold for $\widetilde{G}:$
\begin{wst}
\item[{\rm (i)}] $l\equiv 2\pmod 4;$
\item[{\rm (ii)}] $\eta(\widetilde{C}_l)\equiv 2\pmod 4;$
\item[{\rm (iii)}] $r(T_G)=r([T_G])$.
\end{wst}
\end{prop}
\begin{proof}

For ``sufficiency", we proceed by induction on $|V_{T_G}|$ to show that $\widetilde{G}$ is lower-optimal. If $|V_{T_G}|=1,$ then $G\cong C_l$, by Fact 4, it is straightforward to check that $rk(\widetilde{G})-r(G)=-2.$ Now we assume that $|V_{T_G}|\geqslant 2.$

If $T_G$ is an empty graph, then $\widetilde{G}$ consists of {{the unique mixed cycle $\widetilde{C}_l$}} and isolated vertices. Thus by Facts~2 and 4, $\widetilde{G}$ is lower-optimal.

If $T_G$ is non-empty, together with (iii) and Lemma \ref{lem2}, then we obtain $m(T_G)=m([T_G])$. Note that $T_G$ contains  a pendant vertex, say $x$. Clearly, it is also a pendant vertex of $\widetilde{G}.$ Then, let $y$ be the neighbor of $x.$ By Lemma \ref{lem0.2}, $y$ is not on any cycle of $G.$ Put $\widetilde{G}'=\widetilde{G}-x-y.$ Then by Lemma \ref{lem3}, $r(T_{G'})=r(T_G)-2=r([T_G])-2=r([T_{G'}]).$ Applying the induction hypothesis on $\widetilde{G}'$, we have $\widetilde{G}'$ is lower-optimal. Combining Property~\ref{lem5.02} yields that $\widetilde{G}$ is lower-optimal.

For ``necessity", let $m=m(G)$ be the matching number of $G.$ Combining Lemmas \ref{lem33} and \ref{lem2222}, we obtain $2m-2 \leqslant r(G)\leqslant 2m+1$ and $2m-2 \leqslant rk(\widetilde{G})\leqslant 2m+1.$ Note that $\widetilde{G}$ is lower-optimal, i.e. $rk(\widetilde{G})- r(G)=-2$. Hence, we proceed by considering the following two possible cases.

The first case is $rk(\widetilde{G})=2m-1,$ $r(G)=2m+1.$ In this case, by Lemma \ref{lem33}, we have $l$ is odd. Let $\phi (\widetilde{G},\lambda)=\sum\limits_{i=0}^n a_i\lambda^{n-i}$ be the characteristic polynomial of $\widetilde{G}.$ Then by Lemma \ref{lem7} we obtain $a_{2m}=\sum\limits_{M\in \mathscr{M}} (-1)^{m}=|\mathscr{M}|(-1)^m\neq 0,$ where $\mathscr{M}$ denotes the set of all the maximum matchings of $G.$ Consequently, $rk(\widetilde{G})\geqslant 2m,$ a contradiction.

The rest case is $rk(\widetilde{G})=2m-2,$ $r(G)=2m.$ In view of Lemma \ref{lem2.3} we have $l$ is even, $\eta(\widetilde{C}_l)\equiv l\pmod {4}$ and $m(T_G)=m([T_G])$. Together with Lemma \ref{lem2}, we have $r(T_G)=r([T_G]).$ Hence, (iii) hold. If (i) or (ii) does not hold, then $l\equiv 0\pmod {4}.$ By Lemma \ref{lem33}, $r(G)=2m-2,$ a contradiction.

This completes the proof.
\end{proof}

\setlength{\baselineskip}{17pt}

\begin{prop}\label{lem5.4}
Let $\widetilde{G}$ be obtained from mixed graphs $\widetilde{C}_q$ and $\widetilde{H}$ by connecting them with a mixed edge $\tilde{e}=xy$, where $x\in V_{\widetilde{C}_q}$ and $y\in V_{\widetilde{H}}$. Put $\widetilde{K}:=\widetilde{H}+x.$ If $\widetilde{G}$ is lower-optimal, then
\begin{wst}
\item[{\rm (i)}] for each mixed cycle $\widetilde{C}_l$ of $\widetilde{G},$ $l\equiv 2\pmod {4}$ and $\eta(\widetilde{C}_l)\equiv 2\pmod {4}$;
\item[{\rm (ii)}] both $\widetilde{K}$ and $\widetilde{H}$ are lower-optimal;
\item[{\rm (iii)}] $rk(\widetilde{K})=rk(\widetilde{H})$ and $r(K)=r(H);$
\item[{\rm (iv)}] $rk(\widetilde{G})=q-2+rk(\widetilde{K})$ and $r(G)=q+r(K).$
\end{wst}
\end{prop}
\begin{proof}

(i)\ We begin by induction on $d(G)$ to show our result. If $d(G)=1$, then $\widetilde{G}$ contains a unique mixed cycle and (i) follows directly from Property \ref{lem5.3}. Now suppose that $d(G)\geqslant 2$. Then $\widetilde{H}$ contains at least one mixed cycle. Let $u$ be a vertex lying on some mixed cycle of $\widetilde{H}$. By Property \ref{lem5.1}(ii), $\widetilde{G}_0=\widetilde{G}-u$ is lower-optimal. Since $d(G_0)<d(G)$, by induction each mixed cycle in $\widetilde{G}_0$, including $\widetilde{C}_q$, satisfies (i). By a similar discussion as above, each mixed cycle in $\widetilde{G}-x$ satisfies (i), i.e., all the mixed cycles in $\widetilde{H}$ satisfy (i). That is to say, each mixed cycle in $\widetilde{G}$ satisfies (i).

(ii)\ By (i) we know that $q$ is even. Then let $C_q=xx_2x_3\dots x_qx$. As $\widetilde{G}$ is lower-optimal, by Property \ref{lem5.1} both $\widetilde{G}-x_2$ and $\widetilde{G}-x$ are lower-optimal. Together with Fact \ref{fac1}, $\widetilde{H}$ is also lower-optimal. Let $\widetilde{G}_1:=\widetilde{G}-\{x_2,x_3,x_4\}$. Note that $x_3$ (resp. $x_4$) is the pendant vertex (resp. quasi-pendant vertex) of $\widetilde{G}-\{x_2\}.$ Then in view of Lemma~\ref{lem3}, we have
\begin{eqnarray}\label{eq:5.06}
rk(\widetilde{G}_1)=rk(\widetilde{G}-x_2)-2,\ \ r(G_1)=r(G-x_2)-2.
\end{eqnarray}
As $\widetilde{G}-x_2$ is lower-optimal, we obtain
\begin{eqnarray}\label{eq:5.07}
rk(\widetilde{G}-x_2)-r(G-x_2)=-2d(G-x_2).
\end{eqnarray}
Combining (\ref{eq:5.06}) and (\ref{eq:5.07}) with the fact that $d(G_1)=d(G-x_2)$ yields
$
rk(\widetilde{G}_1)-r(G_1)=-2d(G_1),
$
i.e.,  $\widetilde{G}_1$  is lower-optimal. Repeating such process, after $\frac{q-2}{2}$ steps, the resultant graph is $\widetilde{G}-\{x_2,x_3,\ldots,x_{q}\}=\widetilde{K},$ which is also lower-optimal.

(iii) and (iv)\ By Lemma \ref{lem3} and Property \ref{lem5.1}, one has
\begin{eqnarray}\label{eq:5.7}
&&rk(\widetilde{G})=rk(\widetilde{G}-x_2)=rk(\widetilde{K})+q-2,\ \ \ \ \
r(G)=r(G-x_2)+2=r(K)+q
\end{eqnarray}
and
\begin{eqnarray}\label{eq:5.9}
&& rk(\widetilde{G})=rk(\widetilde{G}-x)=rk(\widetilde{H})+q-2, \ \ \ \ \ r(G)=r(G-x)+2=r(H)+q.
\end{eqnarray}
Together with (\ref{eq:5.7}) and (\ref{eq:5.9}), (iii) and (iv) hold.

This completes the proof.
\end{proof}

\begin{prop}\label{lem5.6}
Let $\widetilde{G}$ be a mixed graph. If $\widetilde{G}$ is lower-optimal, then
\begin{wst}
\item[{\rm (i)}] the cycles (if any) of $G$ are pairwise vertex-disjoint;
\item[{\rm (ii)}] for each mixed cycle $\widetilde{C}_l$ of $\widetilde{G},$ $l\equiv 2\pmod 4$ and $\eta(\widetilde{C}_l)\equiv 2\pmod {4}$;
\item[{\rm (iii)}] $r(G)=r(T_G)+\sum_{O\in \mathcal{C}_G} |V_O|,$ where $\mathcal{C}_G$ denotes the set of cycles of $G;$
\item[{\rm (iv)}] $r(T_G)=r([T_G]).$
\end{wst}
\end{prop}

\begin{proof}
(i)\ By Property \ref{lem5.1}(v), (i) follows immediately.

(ii)-(iv)\ We proceed by induction on ${{|V_G|}}$ to show (ii)-(iv) simultaneously. If $|V_G|=1,$ then (ii)-(iv) hold trivially. Now suppose $|V_G|\geqslant 2.$ By (i) all cycles (if any) of $G$ are pairwise vertex-disjoint. Hence $T_G$ is well defined. If $T_G$ is empty, then $G$ consists of distinct cycles and isolated vertices. Therefore, by Facts \ref{fac1} and \ref{fac2}, (ii)-(iv) hold immediately. If $T_G$ is non-empty, then $T_G$ has a pendant vertex, say $x.$ If $x$ doesn't lie on any cycle of $G,$ then $x$ is also a pendant vertex of $G.$ If $x$ lies on some cycle of $G,$ then this cycle must be a pendant cycle of $G$. We proceed by considering the following two possible cases to show our result.

{\bf{Case 1.}}\ $x$ is a pendant vertex of $G.$ In this case, let $y$ be the unique neighbor of $x$ and put $\widetilde{G}':=\widetilde{G}-x-y.$ Then by Property \ref{lem5.02}, $y$ is not on any cycle of $G$ and $\widetilde{G}'$ is also lower-optimal. Applying induction hypothesis to $\widetilde{G}'$ yields
\begin{wst}
  \item[{\rm (a)}] for each mixed cycle $\widetilde{C}_l$ of $\widetilde{G}',$ $l\equiv 2\pmod 4$ and $\eta(\widetilde{C}_l)\equiv 2\pmod {4}$;
  \item[{\rm (b)}] $r(G')=r(T_{G'})+\sum_{O\in \mathcal{C}_{G'}} (|V_O|);$
  \item[{\rm (c)}] $r(T_{G'})=r([T_{G'}]).$
\end{wst}

Note that each of cycles in $\widetilde{G}$ is also that of $\widetilde{G}'.$ Hence, Assertion (a) implies that each cycle of $\widetilde{G}$ satisfies (ii). Clearly, $x$ is a pendant vertex of $T_G$ (resp., $[T_G]$) whose neighbor is $y$. Remind that $T_{G'}=T_G-x-y$ and $[T_{G'}]=[T_G]-x-y$. Hence, together with Lemma \ref{lem3} and (b), we have
$$
r(G)=r(G')+2=\left[r(T_{G'})+\sum_{O\in \mathcal{C}_{G'}} (|V_O|)\right]+2=r(T_G)+\sum_{O\in \mathcal{C}_G} |V_O|,
$$
i.e., (iii) holds. Combining Lemma \ref{lem3} and (c) yields
$
r(T_G)=r(T_{G'})+2=r([T_{G'}])+2=r([T_G]).
$
Hence, (iv) holds.

{\bf{Case 2.}}\ $x$ lies on a pendant, say $C_q,$ of $G.$ In this case, let $u$ be the unique vertex of $C_q$ of degree 3. Put $\widetilde{H}:=\widetilde{G}-\widetilde{C}_q$ and $\widetilde{K}:=\widetilde{H}+u.$ By Property \ref{lem5.4}, (ii) holds, and both $\widetilde{K}$ and $\widetilde{H}$ are lower-optimal. Furthermore,

\begin{eqnarray}\label{eq:5.11}
r(G)=r(K)+q.
\end{eqnarray}
Applying induction to $\widetilde{K}$ yields
\begin{eqnarray}\label{eq:5.12}
r(K)=r(T_K)+\sum_{O\in \mathcal{C}_{K}} |V_O|.
\end{eqnarray}
Note that $T_G\cong T_K$ and $|V_{C_q}|=q.$ Combining (\ref{eq:5.11}) and (\ref{eq:5.12}) yeilds
$$
r(G)=r(T_G)+\sum_{O\in \mathcal{C}_{G}} |V_O|.
$$
Applying induction to $\widetilde{H}$ yields
\begin{eqnarray}\label{eq:5.13}
r(H)=r(T_H)+\sum_{O\in \mathcal{C}_{H}} |V_O|, \ \ \ r(T_H)=r([T_H]).
\end{eqnarray}
By Property \ref{lem5.4}(iii), one has $r(K)=r(H).$ Combining (\ref{eq:5.12}) with (\ref{eq:5.13}) produces
$
r(T_K)=r(T_H).
$
Together with $[T_G]\cong [T_H],$ we obtain that $r(T_G)=r(T_K)=r(T_H)=r([T_H])=r([T_G]).$

This completes the proof.
\end{proof}

\noindent{\bf Proof of Theorem 1.3}\ \ Now we come back to characterize the sufficient and necessary conditions for the equality holding on the left hand of (\ref{1}).

For ``sufficiency", suppose that $\widetilde{G}$ satisfies all the conditions of (i)-(iii). Then we may apply, repeatedly, $\delta$-transformation for $t$ times to switch $G$ to a crucial subgraph $G_0,$ which is consisted of $d(G)$ disjoint cycles and isolated vertices. Note that each $\delta$-transformation decreases the $H$-rank of a mixed graph by 2, so does the rank of it's underlying graph. Then
\begin{eqnarray}\label{eq:5.8}
rk(\widetilde{G})=2t+rk(\widetilde{G}_0), \ \ \ r(G)=2t+r(G_0).
\end{eqnarray}
By Lemmas \ref{lem1}, \ref{lem10} and \ref{lem12} one has
\begin{eqnarray}\label{eq:5.08}
rk(\widetilde{G}_0)=\sum_{\widetilde{O}\in \mathcal{C}_{\widetilde{G}_0}} rk(\widetilde{O})=\sum_{O\in \mathcal{C}_{G_0}} (|V_O|-2)=\sum_{O\in \mathcal{C}_{G_0}} (r(O)-2)=r(G_0)-2d(G_0).
\end{eqnarray}
Substituting (\ref{eq:5.8}) into (\ref{eq:5.08}) yields $rk(\widetilde{G})-r(G)=-2d(G),$ as desired.

For ``necessity", (i) and (ii) follow directly by Property \ref{lem5.6}. Hence, in order to complete the proof, it suffices to show (iii). We proceed by induction on $|V_{\widetilde{G}}|$ to prove (iii).

If $|V_{\widetilde{G}}|=1,$ then (iii) holds trivially. Suppose that (iii) holds for each lower-optimal mixed graph $\widetilde{G}$ with $|V_{\widetilde{G}}|<n$. Now let $\widetilde{G}$ be an lower-optimal mixed graph of order $n\geqslant 2.$

If $T_G$ is empty, then $\widetilde{G}$ is the disjoint union of $d(G)$ cycles along with some isolated vertices. Then (iii) holds by Facts \ref{fac1} and \ref{fac2}. If $T_G$ is non-empty, then by Property \ref{lem5.6} we have $r(T_G)=r([T_G]).$ Hence, by Lemma \ref{lem5}, $T_G$ contains a pendant vertex, say $x,$ which is also a pendant vertex of $G$. Let $y$ be the neighbor of $x$ and put $\widetilde{G}_1=\widetilde{G}-x-y.$ By Property \ref{lem5.02}, $y$ is not on any mixed cycle of $\widetilde{G}$ and $\widetilde{G}_1$ is lower-optimal. If we apply induction to $\widetilde{G}_1$ we see that a series $\delta$-transformations can switch $\widetilde{G}_1$ to a crucial subgraph $G_0,$ which consists of $d(G_1)$ disjoint cycles together with some isolated vertices. Note that $d(G)=d(G_1).$ Hence, a series of $\delta$-transformation can switch $G$ to the crucial subgraph $G_0,$ which is the disjoint union of $d(G)$ disjoint cycles together with some isolated vertices.

This completes the proof. \qed

\end{document}